\documentclass[11pt,reqno]{amsart}

\usepackage{amsfonts,amssymb}              
\usepackage{bbm}
\usepackage{mathrsfs}
\usepackage{a4wide}

\usepackage{hyperref}

\hypersetup{
  colorlinks   = true, 
  urlcolor     = blue, 
  linkcolor    = blue, 
  citecolor   = red 
}

\newtheorem{neu}{}[section]
\newtheorem{Cor}[neu]{Corollary}
\newtheorem*{Cor*}{Corollary}
\newtheorem{Thm}[neu]{Theorem}
\newtheorem*{Thm*}{Theorem}
\newtheorem*{Observation*}{Observation}

\newtheorem*{Prop*}{Proposition}                
\newtheorem{Lemma}[neu]{Lemma}
\theoremstyle{definition}
\newtheorem*{Rmk*}{Remark}
\newtheorem{Rmk}[neu]{Remark}

\newtheorem*{Ex*}{Example}

\newtheorem*{Qu*}{Question}

\def\1{\:\!}
\def\2{\;\!}
\def\ni{\noindent}

\def\m{\medskip}

\newcommand{\N}{\mathbb{N}}

\newcommand{\Z}{\mathbb{Z}}
\newcommand{\R}{\mathbb{R}}
\newcommand{\C}{\mathbb{C}}

\newcommand{\id}{\mathrm{id}}

\newcommand{\om}{\omega}

\newcommand{\A}{\mathscr{A}}

\newcommand{\D}{\mathbb{D}}

\renewcommand{\L}{\mathscr{L}}
\renewcommand{\H}{\mathrm{H}}
\def\HF{\operatorname{HF}}
\def\SB{\operatorname{SB}}
\def\contr{\operatorname{contr}}

\newcommand{\beq}{\begin{equation}}
\newcommand{\beqn}{\begin{equation}\nonumber}
\newcommand{\eeq}{\end{equation}}

\newcommand{\bea}{\begin{equation}\begin{aligned}}
\newcommand{\bean}{\begin{equation}\begin{aligned}\nonumber}
\newcommand{\eea}{\end{aligned}\end{equation}}

\numberwithin{equation}{section}

\newcommand{\p}{\partial}

\begin{document}
\title[Hamiltonian delay equations]{Hamiltonian delay equations -- examples and a lower bound for the number of periodic solutions}
\author{Peter Albers, Urs Frauenfelder, Felix Schlenk}

\address{Peter Albers\\
 Mathematisches Institut\\
 Ruprecht-Karls-Universit\"at Heidelberg} 
\email{palbers@mathi.uni-heidelberg.de}

\address{Urs Frauenfelder\\
 Mathematisches Institut\\
Universit\"at Augsburg}
\email{urs.frauenfelder@math.uni-augsburg.de}

\address{Felix Schlenk\\
 Institut de Math\'ematiques\\ 
 Universit\'e de Neuch\^atel}
\email{schlenk@unine.ch}


\date{\today}

\begin{abstract}
We describe a variational approach to a notion of Hamiltonian delay equations.
Our delay Hamiltonians are of product form. 
We consider several examples. 
For closed symplectically aspherical symplectic manifolds $(M,\omega)$ we
prove that for generic delay Hamiltonians the number of 1-periodic solutions of 
the Hamiltonian delay equation
is at least the sum of the Betti numbers of~$M$, 
extending the proof of the Arnold conjecture to the case with delay. 
\end{abstract}

\maketitle

\section{Introduction and main results}

An ordinary differential equation (ODE) on $\R^d$ is, in the simplest case, of the form
$$
\dot v(t) = X (v(t))
$$
where $X$ is a vector field on $\R^d$.
A delay differential equation (DDE) on $\R^d$ is, again in the simplest case, of the form
$$
\dot v(t) = X (v(t-\tau))
$$
where $X$ is still a vector field on $\R^d$ and $\tau >0$ is the time delay.
Delay equations therefore model systems in which the instantaneous velocity $\dot v (t)$
depends on the state of the curve~$v$ at a past time. 
There are very many such systems in science and engineering. 
We refer to~\cite{HaVL93} for a foundational text and to~\cite{Erneux09} for a wealth of examples.

A {\it Hamiltonian}\/ differential equation on $\R^{2n}$ is an ODE of the form
\begin{equation} \label{e:Ham}
\dot v (t) = X_H (v(t))
\end{equation}
where the Hamiltonian vector field is of the special form $X_H = i\2 \nabla H$.
Here $H \colon \R^{2n} \to \R$ is a smooth function and $i$ is the usual complex multiplication 
on $\R^{2n} \cong \C^n$.
It is now tempting to define a Hamiltonian delay equation to be a DDE of the form
\begin{equation} \label{e:Liu} 
\dot v (t) = X_H (v(t-\tau))
\end{equation}
with $\tau >0$ and $X_H$ as before.
Such systems where studied by Liu~\cite{Liu12},
who proved the existence of periodic orbits under natural assumptions on~$H$.

In this paper we take a different approach to Hamiltonian delay equations, or at least to periodic
orbits solving what we propose to call a Hamiltonian delay equation.
Our approach is through action functionals.
Let $\L = C^\infty (S^1,\R^{2n})$ be the space of smooth 1-periodic loops in $\R^{2n}$, 
and recall from classical mechanics that the 1-periodic solutions of~\eqref{e:Ham}
are exactly the critical points of the action functional $\A \colon \L \to \R$
given by
$$
\A (v) \,=\, \int_0^1 \bigl[ p(t) \cdot \dot q (t) - H (v(t)) \bigr] \2 dt, \qquad v(t) = (q(t),p(t)).
$$
This fact, that interesting solutions can be seen as the critical points of a functional,
played a key role in the creation of the modern theory of Hamiltonian dynamics and of symplectic topology,
see the outlook at the end of this introduction.
We therefore look at ``delay action functionals''.
If we just take
$$
\A (v) \,=\, \int_0^1 \bigl[ p(t) \cdot \dot q (t) - H (v(t-\tau)) \bigr] \2 dt,
$$
we get nothing new: The critical point equation is again $\dot v (t) = X_H(v(t))$.
However, if we take two Hamiltonian functions $H,K$ on~$\R^{2n}$ and the functional
\begin{equation} \label{e:A}
\A (v) \,=\, \int_0^1 \bigl[ p(t) \cdot \dot q (t) - H(v(t)) \,K (v(t-\tau)) \bigr] \2 dt,
\end{equation}
then the critical point equation is the honest delay equation
\begin{equation} \label{e:Ae}
\dot v (t) \,=\, H(v(t+\tau)) \, X_K(v(t)) + K(v(t-\tau))\, X_H(v(t))   .
\end{equation}
Notice that the functional~\eqref{e:A} makes sense since $v$ is 1-periodic. 
This functional would make no sense on the space of paths $v \colon [0,1] \to \R^{2n}$ 
between two given points, since then $v(t-\tau)$ would not be defined for $t \in [0, \min \{\tau,1\})$.
Also notice that the time shift $+\tau$ in~\eqref{e:Ae} looks like ``into the future'', 
so that equation~\eqref{e:Ae} looks like a forward-backward delay equation.
However, along 1-periodic orbits~$v$ we have $v(t+\tau) = v(t+ (\tau-k))$ for every $k \in \N$,
and so \eqref{e:Ae} is really a delay equation.

In our approach a Hamiltonian delay equation 
is thus a delay equation that can be obtained as critical point equation of
an action functional.

In Sections~\ref{s:2} and~\ref{s:4} we compute the critical point equations of several classes 
of delay action functionals on the loop space of~$\R^{2n}$. 
As a special case we shall obtain in Section~\ref{s:3} 
one instance of the delayed Lotka--Volterra equations. 
In fact, already in his 1928 paper~\cite{Vol28} and in his seminal book \cite{Vol31} from~1931 Volterra was interested in periodic solutions of delay equations, and formulated the famous Lotka--Volterra equations with and without delay. 
In Section~\ref{s:int} we give a first integral along periodic orbits for certain
Hamiltonian delay equations.

\m
\ni
{\bf Extension to manifolds.}
A symplectic manifold is a manifold~$M$ together with 
a non-degenerate closed 2-form~$\omega$ on~$M$.
We assume throughout that $(M,\omega)$ is symplectically aspherical, 
meaning that the cohomology class~$[\omega]$ vanishes on the image 
of~$\pi_2(M)$ in~$H_2(M)$ under the natural map forgetting the base point.
Examples of symplectically aspherical manifolds are exact symplectic manifolds, 
for which $\omega = d \lambda$ for a 1-form~$\lambda$, 
like $\R^{2n}$ or cotangent bundles with their canonical symplectic form, 
and closed examples are closed orientable surfaces of positive genus and their products, 
like tori~$T^{2n}$.

Recall that on $\R^{2n}$ one possible definition of a Hamiltonian delay equation is
\begin{equation} \label{e:Hamtau}
\dot v(t) = X_H(v(t-\tau)) \;.
\end{equation}
On a general symplectic manifold~$M$, however, this concept does not make sense, 
simply because $\dot v(t) \in T_{v(t)} M$ and $X_H(v(t-\tau)) \in T_{v(t-\tau)}M$ 
reside in different tangent spaces. 
On the other hand, our approach through action functionals readily extends to manifolds:
The Hamiltonian vector field of a smooth function $H \colon M \to \R$ is defined by
$\omega (X_H, \cdot ) = -dH$, 
and the contractible 1-periodic solutions of~$X_H$ are exactly the critical points
of the action functional~$\A$ defined on the component of contractible loops~$\L_{\contr}$ 
by
\begin{equation} \label{e:AHo}
\A (v) \,=\, \int_{\D} \bar v^*\omega  - \int_{0}^1 H (v(t)) \2 dt
\end{equation}
where $\bar v \colon \D \to M$ is a smooth map on the unit disc such that $\bar v(e^{2\pi i t}) = v(t)$.
The value of $\int_{\D} \bar v^*\omega$ does not depend on our choice of the filling disc~$\bar v$
in view of the asphericity assumption on~$(M,\omega)$.
In the special case that $\omega = d\lambda$ is exact, the action functional is
defined on the full loop space~$\L$ by
\begin{equation} \label{e:AH}
\A (v) \,=\, \int_0^1 \bigl[ \lambda (\dot v) - H (v(t)) \bigr] \2 dt .
\end{equation}
Taking $\lambda = \sum_{j=1}^n p_j\1 dq_j$ on $\R^{2n}$ we recover the case
described before.

Replacing the Hamiltonian term $H(v(t))$ in~\eqref{e:AHo} by the delay term  
$H(v(t)) \,K (v(t-\tau))$, we get as critical point equation the delay equation~\eqref{e:Ae} on~$M$,
and for any of the terms described in Sections~\ref{s:2} and~\ref{s:4}
we get other Hamiltonian delay equations on~$M$.
Thus, if we start from a delay action functional and compute the critical point equation, 
then an ``accident'' as for~\eqref{e:Hamtau} cannot happen, 
and we always get a meaningful equation.

The search for periodic orbits is a central theme in Hamiltonian dynamics. 
The modern tools of symplectic geometry, in particular Floer homology, 
imply lower bounds for the number of 1-periodic orbits of Hamiltonian systems
that are 1-periodic in time. 
In particular, the Arnold conjectures are proven on all symplectically aspherical manifolds, 
see~\cite[\S 12]{McSal12} and~\cite{AbSch19} for a survey on Floer homology. 
Do these lower bounds also hold for Hamiltonian delay equations?
The main result of this paper is that this is so for the special class of Hamiltonian delay 
equations in product form. 
Given a closed manifold $M$ of dimension~$d$ denote by
$$
\SB (M) \,=\, \sum_{j=0}^d b_j (M;\Z_2)
$$
the total dimension of the $\Z_2$-homology of~$M$.

\begin{Thm} \label{t:floer}
Assume that $(M,\omega)$ is a closed symplectically aspherical symplectic manifold.
Fix a time delay $\tau >0$, and let $H,K \colon M \times S^1 \to \R$ be two time-dependent Hamiltonian functions such that
the Hessian of the action functional $\A \colon \L_{\contr} \to \R$,
$$
\A (v) \,=\, \int_{\D} \bar v^*\omega 
          - \int_0^1 H \big(v(t),t\big) K \big(v(t-\tau), t-\tau \big) \2 dt ,
$$
at its critical points $v$ has trivial kernel. 
Then $\A$ has at least $\SB (M)$ critical points. 
In other words, there are at least $\SB (M)$ contractible 1-periodic orbits in~$M$ 
that solve the Hamiltonian delay equation
$$\dot v (t) \,=\, H(v(t+\tau)) \, X_K(v(t)) + K(v(t-\tau))\, X_H(v(t)) .
$$
\end{Thm}

Note that taking $K$ to be a constant function we obtain the classical Arnold conjecture
as a special case.
It will be clear from the proof in~\S \ref{s:floer} that Theorem~\ref{t:floer} 
also holds for all the examples of Hamiltonian delay equations given in Section~\ref{s:4}.

\m \ni
{\bf Related works.}
General properties of delay action functionals were studied for instance 
in Chapter~VI of~\cite{Els64} and in~\cite{Sab69}.
The problem when a delay equation on~$\R^d$ is the critical point equation of a functional 
is analyzed in~\cite{KPS07}.
A Hamiltonian formalism for certain non-local PDEs on~$\R^{2n}$,
that is also based on non-local action functionals, 
was recently proposed in~\cite{BaSch17}.

General DDEs can be readily defined on manifolds, see Section~12.1 of~\cite{HaVL93} and~\cite{Oli69}.
In contrast, there seems to be no concept of a Hamiltonian delay equation on manifolds.
Our approach at least provides a natural notion of such an equation 
and a tool for finding periodic solutions.

Results on periodic orbits of DDEs are very scarce, even in~$\R^d$,
see the end of the introduction to~\cite{AFS2}, and \cite{BCF09} and the references therein.
The only previous multiplicity results for Hamiltonian DDEs
are for a certain class of equations of the form~\eqref{e:Liu} in~$\R^{2n}$, \cite{Liu12},
and for a special class of Hamiltonian delay equations with piece-wise smooth orbits
on symplectically aspherical manifolds, \cite{AFS2},
where the Arnold conjectures were derived from classical Lagrangian Floer homology
by an iterated graph construction.

\m \ni
{\bf Outlook: A calculus of variations for Hamiltonian DDEs.}
In the rest of this introduction we further comment on why we believe that our approach
to Hamiltonian delay equations is promising.

While the theory of DDEs is meanwhile quite rich~\cite{HaVL93}, 
it is nonetheless much less developed than the theory of ODEs. 
One reason is that for DDEs there is no local flow on the given manifold,
and so the whole theory is less geometric and more cumbersome. 

In sharp contrast to general ODEs, Hamiltonian ODEs can be studied by variational methods, 
thanks to the action functional.
For one thing, the action functional 
(even though neither bounded from below nor above, and strongly indefinite) 
can be used to do critical point theory on the loop space, 
as was first demonstrated by Rabinowitz~\cite{Rab78}.
At least as important, one can see
symplectic topology as the geometry of the action functional.
(This is almost the title and exactly the content of Viterbo's paper~\cite{Vi92}, 
and also in the book~\cite{HoZe94} the action functional is the main tool.)
For instance, a selection of critical values of the action functional by min-max leads
to numerical invariants of Hamiltonian systems and symplectic manifolds, that have many applications.
The climactic impact of the action functional into symplectic dynamics and topology, however, 
is Floer homology, which is Morse theory for the action functional on the loop space.

Now, incorporating delays into an action functional, we can try to extend all these constructions
to delay action functionals, 
and thereby create a calculus of variations for Hamiltonian delay equations, that should have many applications
at least to questions on periodic orbits of such equations.
We have done this in the present article for special Hamiltonian delay equations, and
our proof of Theorem~\ref{t:floer} applies to all Hamiltonian delay equations 
in which the delay~$\tau$ is constant or just depends on time, 
but not on the loop~$v$.
There are, however, many interesting delay equations in which the delay depends 
on the loop. For instance, in the equation considered by Carl Neumann in~1868 
in his derivation of Weber's force law in electrodynamics, 
the delay is of the form $\tau (v) = - \frac{|q(t)|}{c}$ for a loop $v(t) =(q(t),p(t))$ in~$T^*\R^2$, 
see~\cite{Fra19, FraWeb18}.
It is an interesting problem to understand whether for such delay equations
a Floer homology can still be constructed.
First steps in the construction of a general delay Floer homology 
were taken in~\cite{AlFr13} and~\cite{AFS1}.

\subsubsection*{Acknowledgment}
PA is supported by DFG CRC/TRR 191 and under Germany’s Excellence Strategy EXC-2181/1 - 390900948 (the Heidelberg STRUCTURES Excellence Cluster).
UF is supported by DFG FR/2637/2-1, and 
FS is supported by the SNF grant 200020-144432/1.

\section{Delay equations from sums and products of Hamiltonian functions} \label{s:2}

Let $(M,\omega)$ be a symplectically aspherical manifold, i.e., $[\omega] |_{\pi_2(M)} =0$.
We choose $2N+1$ autonomous Hamiltonian functions
\begin{equation*}
F,\, H_i,\, K_i \colon M \to \R, \quad i=1,\ldots,N .
\end{equation*}
Let $\L \equiv \L(M) := C^\infty(S^1,M)$ be the free loop space,
and denote by $\L_{\contr}$ the component of contractible loops.
Now consider the ``action functional'' $\A  \colon \L_{\contr} \to \R$,
\begin{equation} \label{eqn:action_fctl_with_sum}
\A(v) \,=\, \int_{\D} \bar v^*\omega - \int_0^1 F \big( v(t) \big)dt 
          - \sum_{i=1}^N \int_0^1 H_i \big(v(t)\big) K_i \big(v(t-\tau_i) \big) dt
\end{equation}
where $\bar v \colon \D \to M$ is a smooth map on the unit disc such that $\bar v(e^{2\pi i t}) = v(t)$
and where $\tau_1, \dots, \tau_N \geq 0$ are $N$ time delays. 
The value of $\int_{\D} \bar v^*\omega$ does not depend on our choice of the filling disc~$\bar v$
in view of the asphericity assumption on~$(M,\omega)$.
If $\omega = d \lambda$ is exact, we can define $\A$ on all of~$\L$,
by replacing $\int_{\D} \bar v^*\omega$ by $\int_{S^1} v^*\lambda = \int_0^1 \lambda (\dot v (t)) \2 dt$.

To find the critical point equation of~$\A$ we fix $v \in \L_{\contr}$ and $\hat v \in T_v \L_{\contr}$
and compute
\begin{eqnarray*}
d \A(v) \, \hat v &=& \int_0^1 \om \bigl( \hat v(t), \dot v(t) \bigr) dt
       - \int_0^1 dF \bigl( v(t) \bigr) \bigl[ \hat v(t) \bigr] dt \\
&& \;\; -\sum_{i=1}^N \int_0^1 H_i \bigl( v(t) \bigr) \cdot dK_i \bigl( v(t-\tau_i) \bigr) 
                            \bigl[ \hat v(t-\tau_i) \bigr] dt \\
&&\;\;-\sum_{i=1}^N \int_0^1 K_i \bigl( v(t-\tau_i) \bigr) \cdot dH_i \bigl( v(t) \bigr) 
                            \bigl[ \hat v(t) \bigr] dt \,.
\end{eqnarray*}
Since $v$ and $\hat v$ are 1-periodic, 
$$
\int_0^1 H_i \bigl( v(t) \bigr) \cdot dK_i \bigl( v(t-\tau_i) \bigr) \bigl[ \hat v(t-\tau_i) \bigr] dt
\,=\,
\int_0^1 H_i \bigl( v(t+\tau_i) \bigr) \cdot dK_i \bigl( v(t) \bigr) \bigl[ \hat v(t) \bigr] dt.
$$
Using also the definition $\om (X_F,\cdot)=-dF$ 
of the Hamiltonian vector field $X_F$, etc., we find
\begin{eqnarray*}
d \A(v) \, \hat v 
&=& \int_0^1 \om \bigl( \hat v(t), \dot v(t) \bigr) dt 
                  - \int_0^1 \om \bigl( \hat v(t), X_F(v(t)) \bigr) dt \\
&&\;\;-\sum_{i=1}^N \int_0^1 H_i \bigl( v(t+\tau_i) \big) \cdot 
                     \om \bigl( \hat v(t), X_{K_i} (v(t)) \bigr) dt\\
&&\;\;-\sum_{i=1}^N \int_0^1 K_i \bigl( v(t-\tau_i) \bigr) 
                    \cdot \om \bigl( \hat v(t), X_{H_i} (v(t)) \bigr) dt .
\end{eqnarray*}
The critical point equation is therefore
\begin{equation*}
\dot v(t) \,=\, X_F(v(t)) + 
\sum_{i=1}^N \Bigl[ H_i (v(t+\tau_i)) \, X_{K_i} (v(t)) + K_i (v(t-\tau_i)) \, X_{H_i}(v(t)) \Bigr] .
\end{equation*}
We have proved the following lemma.
\begin{Lemma} \label{le:crit}
The critical points of $\A$ satisfy the Hamiltonian delay equation
\begin{equation} \label{eqn:Ham_delay_eqn}
\dot v(t) \,=\, X_F(v(t)) + 
   \sum_{i=1}^N \Bigl[ H_i(v(t+\tau_i)) \, X_{K_i} (v(t)) + K_i (v(t-\tau_i)) \, X_{H_i}(v(t)) \Bigr] .
\end{equation}
\end{Lemma}

Using that $v(t+1)=v(t)$ we obtain

\begin{Cor}
For the single time delay $\tau=\tfrac12$ the Hamiltonian delay equation becomes
\begin{equation*} 
\dot v(t) \,=\, X_F (v(t)) + 
      \sum_{i=1}^N \Bigl[ H_i (v(t-\tfrac12)) \, X_{K_i}(v(t)) + 
	              K_i(v(t-\tfrac12)) \, X_{H_i}(v(t)) \Bigr] .
\end{equation*}
\end{Cor}

\section{The Lotka--Volterra equations, with and without delay} \label{s:3}

In this section we extend the work of Fernandes--Oliva~\cite{FeOl95} 
to positive delays. 
Fix a skew-symmetric $N \times N$-matrix $A=(a_{ij})$, i.e.~$a_{ji}=-a_{ij}$,
and $N$ real numbers~$b_i$. 
Take $M = \R^{2N}$ with its usual exact symplectic form $\om = \sum_i dq_i \wedge dp_i$,
and set
\begin{equation*}
F(q,p) := \sum_{i=1}^N b_i q_i, \qquad
H_i(q,p) := -e^{p_i}, \qquad
K_i(q,p) := e^{\frac12 \sum_{j=1}^N a_{ij}q_j} .
\end{equation*}
The Hamiltonian vector fields are
\begin{equation*}
X_F = \sum_{i=1}^N b_i \frac{\p}{\p p_i}, \qquad
X_{H_i} = e^{p_i} \frac{\p}{\p q_i}, \qquad
X_{K_i} = \tfrac 12 \sum_{k=1}^N a_{ik} \, e^{\frac12 \sum_{j=1}^N a_{ij}q_j} \frac{\p}{\p p_k}.
\end{equation*}
Fix $\tau \geq 0$.
For $v=(q,p) \in \L(\R^{2N})$ the Hamiltonian (delay) equation~\eqref{eqn:Ham_delay_eqn} 
with equal time delays $\tau_i = \tau$ becomes
\begin{equation*} 
\dot v(t) \,=\, 
\sum_{i=1}^N b_i \frac{\p}{\p p_i} + 
\sum_{i=1}^N \Biggl[-e^{p_i(t+\tau)} \, \tfrac12 \sum_{k=1}^N a_{ik} \, 
                e^{\frac12 \sum_{j=1}^N a_{ij} q_j(t)} \frac{\p}{\p p_k} + 
		    e^{\frac12 \sum_{j=1}^N a_{ij} q_j(t-\tau)} \, e^{p_i(t)} \frac{\p}{\p q_i} \Biggr]\;.
\end{equation*}
In other words,
\begin{eqnarray*}
\dot q_i(t) &=& e^{p_i(t) + \frac12 \sum_{j=1}^N a_{ij} q_j(t-\tau)} \\ [1ex]
\dot p_i(t) &=& b_i- \sum_{l=1}^N e^{p_l(t+\tau)} \, \tfrac12 \, a_{li} \, e^{\frac12\sum_{j=1}^Na_{lj}q_j(t)} \\
&=&b_i- \tfrac12\sum_{l=1}^N a_{li} \, \dot q_l(t+\tau) \\
&=&b_i+ \tfrac12\sum_{l=1}^N a_{il} \, \dot q_l(t+\tau) 
\end{eqnarray*}
where in the last equation we have used that $A$ is skew-symmetric. 
Using these two equations we compute
\begin{eqnarray*}
\ddot q_i(t) &=& 
    \Bigl(\dot p_i(t) + \tfrac12 \sum_{j=1}^N a_{ij} \, \dot q_j(t-\tau) \Bigr) \,\dot q_i(t) \\
&=& \Bigl(b_i+ \tfrac12 \sum_{l=1}^N a_{il} \, \dot q_l(t+\tau) + 
                          \tfrac12 \sum_{j=1}^N a_{ij} \, \dot q_j(t-\tau) \Bigr) \,\dot q_i(t)\;.
\end{eqnarray*}
Now observe that the right hand side only depends on the $\dot q_j$, 
but not on the $p_j$.
Setting $x_i(t) := \dot q_i(t)$ we thus obtain the first order delay system
\begin{equation} \label{e:1oDDE}
\dot x_i(t) \,=\, b_i \, x_i(t) + \tfrac12 \sum_{j=1}^N a_{ij} \, 
           x_i(t) \, x_j(t+\tau) + \tfrac12 \sum_{j=1}^N a_{ij} \, x_i(t) \, x_j(t-\tau)\;.
\end{equation}

\noindent
{\bf Case $\tau =0$.} 
Then \eqref{e:1oDDE} becomes
\begin{equation*}
\dot x_i(t) \,=\, b_i \, x_i(t) + \sum_{j=1}^N a_{ij}\, x_i(t)\, x_j(t)
\end{equation*}
with skew-symmetric $A=(a_{ij})$. 
This is one instance of the Lotka--Volterra equations without delay. 
These equations were proposed by Lotka~\cite{Lot10}
in his studies of chemical reactions, and independently by Volterra~\cite{Vol26} 
in his studies of predator-prey dynamics.

\medskip
\noindent
{\bf Case $\tau = \frac 12$.} 
Then the equations~\eqref{e:1oDDE} for 1-periodic orbits become 
\begin{equation*}
\dot x_i(t) \,=\,
           b_i \, x_i(t) + \sum_{j=1}^N a_{ij} \, x_i(t) \, x_j(t-\tfrac 12) \,.
\end{equation*}
These Hamiltonian delay equations already appeared in Chapter~4 of Volterra's book~\cite{Vol31}.


\section{More examples} \label{s:4}

In this section we give four rather special classes of Hamiltonian delay equations,
two involving integrals.  
The reader may invent his own examples.
We assume throughout this section that $(M,\omega)$ is symplectically aspherical.

\subsection{Integrals of products of Hamiltonian functions}

In \eqref{eqn:action_fctl_with_sum} we may replace the sum by an integral and choose a double time-dependence: Consider functions $H,K \colon M \times S^1 \times S^1 \to \R$, 
which we write as $H_{t,\tau}(x)$ and $K_{t,\tau}(x)$ for $x \in M$ and $t,\tau \in S^1$. 
Then define $\A \colon \L_{\contr} \to \R$ by
\begin{equation*} 
\A (v) \,=\, \int_{\D} \bar v^*\omega - \int_0^1 \int_0^1 H_{t,\tau} 
      \bigl( v(t) \bigr) \, K_{t,\tau} \bigl( v(t-\tau) \bigr) \,d\tau dt\;.
\end{equation*}
For $v \in \L_{\contr}$ and $\hat v \in T_v \L_{\contr}$ we compute
\begin{eqnarray*}
d \A(v)\, \hat v &=& \int_0^1 \om \bigl(\hat v(t), \dot v(t) \bigr) dt \\
&&\;\; -\int_0^1 \bigg[ \int_0^1 H_{t,\tau} \bigl( v(t) \bigr) \cdot 
                 dK_{t,\tau} \bigl( v(t-\tau) \bigr) \bigl[ \hat v(t-\tau)\bigr] d\tau \\
&&\;\, \phantom{\int_0^1 \bigg[}  - \int_0^1 K_{t,\tau} \bigl( v(t-\tau) \bigr) \cdot 
                 dH_{t,\tau} \bigl( v(t) \bigr) \bigl[ \hat v(t) \bigr] d\tau \bigg] dt \,.
\end{eqnarray*}
Since $v, \hat v$ are 1-periodic and also $H,K$ are periodic in~$t$, 
$$
\int_0^1 H_{t,\tau} \bigl( v(t) \bigr) \cdot 
                 dK_{t,\tau} \bigl( v(t-\tau) \bigr) \bigl[ \hat v(t-\tau)\bigr] d\tau
\,=\,
\int_0^1 H_{t+\tau,\tau} \bigl( v(t+\tau) \bigr) \cdot 
                 dK_{t+\tau,\tau} \bigl( v(t) \bigr) \bigl[ \hat v(t)\bigr] d\tau .		
$$
Therefore, 
\begin{eqnarray*}
d \A(v)\, \hat v &=& \int_0^1 \om \bigl( \hat v(t), \dot v(t) \bigr) dt \\
&& \;\;-\int_0^1\bigg[ \int_0^1 \om \Bigl( \hat v(t), H_{t+\tau,\tau} \bigl(v(t+\tau)\bigr) 
     \cdot X_{K_{t+\tau,\tau}} \bigl( v(t) \bigr) \Bigr) d\tau \\
&& \;\, \phantom{\int_0^1 \bigg[}  -\int_0^1 \om \Bigl( \hat v(t), K_{t,\tau} \bigl(v(t-\tau) \bigr) \cdot 
      X_{H_{t,\tau}} \bigl(v(t)\bigr) \Bigr) d \tau \bigg] dt \,.
\end{eqnarray*}
Hence the critical points of $\A$ are the solutions of the Hamiltonian delay equation
\begin{equation} \label{eqn:Hamiltonian_delay_from_a_product}
\dot v(t) \,=\, \int_0^1 \Bigl[ H_{t+\tau,\tau} \bigl( v(t+\tau) \bigr) \cdot 
            X_{K_{t+\tau,\tau}} \bigl(v(t)\bigr) + K_{t,\tau} \bigl(v(t-\tau)\bigr) \cdot 
		X_{H_{t,\tau}} \bigl(v(t)\bigr)\Bigr] \, d \tau\;.
\end{equation}
In the special case that $H_{t,\tau}$ and $K_{t,\tau}$ are autonomous, 
equation~\eqref{eqn:Hamiltonian_delay_from_a_product} simplifies to
\begin{equation*}
\dot v(t) \,=\, \int_0^1 H \bigl( v(t+\tau) \bigr) \,d\tau \cdot 
       X_K \bigl(v(t)\bigr) + \int_0^1 K \bigl( v(t-\tau) \bigr) \,d\tau \cdot 
	 X_H \bigl(v(t)\bigr) \;.
\end{equation*}
If we define the functions $\overline H, \overline K \colon \L \to \R$ by
\begin{equation*}
\overline H(v) := \int_0^1 H \bigl( v(t) \bigr)\, dt \quad \text{and} \quad 
\overline K(v) := \int_0^1 K \bigl(v(t)\bigr)\, dt ,
\end{equation*}
the above equation becomes
\begin{equation*} 
\dot v(t) \,=\, \overline H(v) \, X_K(v(t)) + \overline K (v) \, X_H(v(t))\;.
\end{equation*}
Specializing further to $H=K$ we obtain
\begin{equation} \label{e:2oH}
\dot v(t) \,=\, 2 \2 \overline H(v) \,X_H(v(t))\;.
\end{equation}
In this special case, preservation of energy implies that $t\mapsto H(v(t))$ is constant along solutions, and thus we may write~\eqref{e:2oH} as a usual Hamiltonian equation:
\begin{equation} \label{eqn:Hamiltonian_delay_autonomous}
\dot v(t) \,=\, 2 \2  H(v) \,X_H(v(t)) \,=\, X_{H^2} (v(t))\;.
\end{equation}

\begin{Rmk}
Of course, this autonomous differential equation can be studied by Floer theory, hence there are 
(in the Morse--Bott sense) multiplicity results (in terms of cup-length or Betti numbers)
for periodic solutions in a certain range of Conley--Zehnder indices. 
Clearly, \eqref{eqn:Hamiltonian_delay_autonomous} has many solutions, namely critical points of~$H$.
However, unless $H$ is $C^2$-small at all critical points, the Morse indices of critical points cannot all agree with their Conley--Zehnder indices, and so Floer theory implies the existence of additional non-constant solutions to~\eqref{eqn:Hamiltonian_delay_autonomous}. 
As we shall see in \S\ref{s:floer}, also equation~\eqref{eqn:Hamiltonian_delay_from_a_product} 
admits a Floer theory. 
Typically, \eqref{eqn:Hamiltonian_delay_from_a_product} has no constant solutions at all, 
even when $K=1$ and $H$ is independent of~$\tau$.
Hence Floer theory implies the existence of many interesting periodic solutions 
to equation~\eqref{eqn:Hamiltonian_delay_from_a_product}.
For instance, the proof of Theorem~\ref{t:floer} in \S\ref{s:floer} shows that if 
$M$ is closed, then 
for generic $H,K$, equation~\eqref{eqn:Hamiltonian_delay_from_a_product}
has at least $\SB (M)$ 1-periodic solutions. 
\end{Rmk}

\subsection{Products of more than two functions}

Generalizing \eqref{eqn:action_fctl_with_sum} in another direction, 
we take $N \geq 3$ functions $H_1, \dots, H_N \colon M \times S^1 \to \R$
and non-negative time delays $\tau_1, \dots, \tau_N$. Then the critical points of 
the functional
\begin{equation*} 
v \,\mapsto\, \int_{\D} \bar v^*\omega - 
    \int_0^1  \prod_{j=1}^N H_j \bigl(v(t-\tau_j)\bigr)  \, dt
\end{equation*}
are the 1-periodic solutions of the Hamiltonian delay equation
$$
\dot v(t) \,=\, \sum_{j=1}^N 
\left[ 
\prod_{\substack{k=1 \\ k \neq j}}^N  H_k 
\bigl( v (t-\tau_k+\tau_j) \bigr) 
\right] 
X_{H_j}(v(t)) .
$$

\subsection{Exponentials of Hamiltonian functions}

We consider yet another incarnation of a Hamiltonian delay equation.
Take
\begin{equation*} 
\begin{aligned}
\A \colon \L_{\contr} & \to \R \\
v&\mapsto \int_{\D} \bar v^*\omega - 
      \int_0^1 \exp \left[ \int_0^1 H_\tau \bigl( v(t-\tau) \bigr) d\tau \right] dt
\end{aligned}
\end{equation*}
where $H \colon M \times S^1 \to \R$ is given. 
We compute
\begin{eqnarray*}
d\A(v) \, \hat v &=& \int_0^1 \om \bigl( \hat v(t), \dot v(t) \bigr) dt \\
&& \;\;-\int_0^1 \int_0^1 \exp \left[\int_0^1 H_\tau \bigl( v(t-\tau) \bigr) d\tau \right]
                 dH_\sigma \bigl( v(t-\sigma) \bigr) \hat v(t-\sigma) \,dt \2 d\sigma .
\end{eqnarray*}
Substituting $t$ by $t+\sigma$ and changing the order of integration the second summand becomes
\begin{eqnarray*}
&& -\int_0^1 \int_0^1 \exp \left[\int_0^1 H_\tau \bigl( v(t+\sigma-\tau) \bigr) d\tau \right]
                          dH_\sigma \bigl( v(t) \bigr) \2 \hat v(t) \, d\sigma \2 dt \\
&=& 				
-\int_0^1 \om \left( \hat v(t), \int_0^1 \exp\left[\int_0^1 H_\tau (v(t+\sigma-\tau)) d\tau \right]                      X_{H_\sigma} \bigl(v(t)\bigr) \, d\sigma \right) dt .
\end{eqnarray*}
The critical point equation is therefore
\begin{equation*}
\dot v(t) \,=\, 
\int_0^1 \exp\left[\int_0^1 H_\tau \bigl( v(t+\sigma-\tau) \bigr) d\tau \right]
                                 X_{H_\sigma} \bigl(v(t)\bigr) \, d\sigma \,.
\end{equation*}

\subsection{Several inputs}

We now consider a function $H \colon M \times M \to\R$ on the symplectic manifold 
$(M\times M,\om \oplus\om)$, 
and denote by $d_1H(x,y) \colon T_xM \to\R$ the derivative 
of~$H$ with respect to the first variable and correspondingly by $X^1_H(x,y)$ 
the Hamiltonian vector field of~$H$ with respect to the first variable:
\begin{equation*}
d_1H(x,y)\xi \,=\, -\om_x \left( X^1_H(x,y),\xi \right) \quad \forall\, \xi \in T_xM\;.
\end{equation*}
Further, we consider the action functional $\A \colon \L_{\contr}(M) \to \R$,
\begin{equation*}
\A (v) \,=\, \int_{\D}\bar v^*\omega - \int_0^1 H \bigl( v(t),v(t+\tau) \bigr) dt\;.
\end{equation*}
Concrete examples for the function~$H$ come
for instance from interaction potentials (as in the $2$-body problem)
or from vortex equations with delay. 
For $v \in \L_{\contr} (M)$ and $\hat v \in T_v \L_{\contr}(M)$ we compute
\begin{eqnarray*}
d\A(v) \,\dot v &=& \int_0^1 \om_{v(t)} \bigl( \hat v(t), \dot v(t) \bigr) \,dt \\
&& \;\;-\int_0^1 \Bigl( d_1H \bigl(v(t),v(t+\tau)\bigr)\, \hat v(t) + 
                d_2H \bigl( v(t),v(t+\tau) \bigr)\, \hat v(t+\tau) \Bigr) \,dt 
\end{eqnarray*}
The second summand is equal to
\begin{eqnarray*}
&&
-\int_0^1 \Bigl( d_1H \bigl(v(t),v(t+\tau)\bigr)\, \hat v(t) + 
                d_2H \bigl( v(t-\tau),v(t) \bigr)\, \hat v(t) \Bigr) \, dt \\
&=&
-\int_0^1 \left[ \om_{v(t)} \Bigl( \hat v(t), X^1_H \bigl( v(t), v(t+\tau) \bigr) \Bigr) + 
               \om_{v(t)} \Bigl( \hat v(t), X^2_H \bigl( v(t-\tau),v(t) \bigr) \Bigr) \right] \,dt .
\end{eqnarray*}
The critical point equation for $\A$ is therefore
\begin{equation} \label{e:critsev}
\dot v(t) \,=\, X^1_H \bigl( v(t),v(t+\tau) \bigr) + X^2_H \bigl( v(t-\tau),v(t) \bigr)\;.
\end{equation}
We point out that indeed
\begin{equation*}
X^1_H \bigl(v(t),v(t+\tau)\bigr), \; X^2_H \bigl( v(t-\tau),v(t) \bigr) \in T_{v(t)}M
\end{equation*}
so that \eqref{e:critsev} makes sense.

\section{A first integral for periodic solutions}
\label{s:int}

For autonomous Hamiltonian systems, the Hamiltonian function is a first integral, 
that is, constant along the solution curves. 
For Hamiltonian delay equations, there is no such function in general.
For instance, consider the Harmonic oscillator $H(z) = \frac 12 \|z\|^2$
on~$\R^2$. Then $X_H (z) = i z$.
We look for solutions of the Hamiltonian delay equation
$$
\dot z (t) \,=\, X_H (z(t-1)) = i \, z(t-1) , \quad t \in \R .
$$ 
Making the Ansatz $z(t) = e^{t v}z(0)$ for $v \in \C$,
we find that $v$ must solve the Lambert equation
$v e^{v} = i$.
By basic properties of the Lambert $W$-function, 
there is a sequence of solutions $v_n$ in~$\C$
with $\mbox{Re}\, v_n \to -\infty$.
Take one such solution~$v$.
Then the time-$t$-map $\varphi_{v}^t(z) = e^{tv}z$
takes the disc $B^2(R)$ to the disc $B^2 ( R \, e^{t \,\mbox{\tiny Re}\;\! v})$ 
that shrinks to the origin as $t \to \infty$. 
The only possible first integrals are therefore constant and thus useless. 

\m
For some of the autonomous Hamiltonian delay equations studied in this paper 
we at least have a first integral along {\it periodic}\/ solutions:
Let $(M,\omega)$ be a symplectically aspherical manifold. 
In equation~\eqref{eqn:Ham_delay_eqn} we take $F=0$ and equal Hamiltonians $H_i = H$ and $K_i = K$
as well as equidistant time delays $\tau_i = \frac iN$, $i=0,\dots,N-1$.
We thus look at periodic solutions $v \in C^\infty(S^1,M)$ of the Hamiltonian delay equation
\begin{equation} \label{hamdel}
\dot v(t) \,=\, 
\Bigg( \sum_{i=0}^{N-1} H\big(v\big(t+\tfrac{i}{N}\big)\big)\Bigg)X_K(v(t))+\Bigg(\sum_{i=0}^{N-1}K\big(v\big(t+\tfrac{i}{N}))\Bigg)X_H(v(t)).
\end{equation}

\begin{Lemma}
The quantity
$$
I_v(t) \,=\, 
\sum_{j,k=0}^{N-1} H \bigl(v\big(t+\tfrac{j}{N}\big)\bigr) \, K\bigl(v\big(t+\tfrac{k}{N}\big)\bigr)
$$
is independent of $t \in S^1$, i.e., $I_v$ is a first integral of the Hamiltonian delay equation~\eqref{hamdel}.
\end{Lemma}

For the proof we show that the derivative of $I_v$ vanishes:
Using the Leibniz rule and plugging in~\eqref{hamdel} we find
\begin{eqnarray*}
\frac{dI_v(t)}{dt} &=& 
\sum_{j,k=0}^{N-1} H\big(v\big(t+\tfrac{j}{N}\big)\big) \,dK\big(v\big(t+\tfrac{k}{N}\big)\big)
                  \, \dot v \big(t+\tfrac{k}{N}\big) + \\
& &\sum_{j,k=0}^{N-1} K\big(v\big(t+\tfrac{k}{N}\big)\big) \, dH\big(v\big(t+\tfrac{j}{N}\big)\big)
                  \, \dot v\big(t+\tfrac{j}{N}\big)\\
&=&\sum_{i,j,k=0}^{N-1} H\big(v\big(t+\tfrac{j}{N}\big)\big) \, K\big(v\big(t+\tfrac{k+i}{N}\big)\big) \, 
                 dK\big(v\big(t+\tfrac{k}{N}\big)\big) \, X_H\big(v\big(t+\tfrac{k}{N}\big)\big)+\\
& &\sum_{i,j,k=0}^{N-1} K\big(v\big(t+\tfrac{k}{N}\big)\big) \, H\big(v\big(t+\tfrac{j+i}{N}\big)\big)\, 
                       dH\big(v\big(t+\tfrac{j}{N}\big)\big) \, X_K \big(v\big(t+\tfrac{j}{N}\big)\big) 
\end{eqnarray*}
In the first sum we switch the indices $i$ and $k$, 
and in the second sum we switch the indices $i$ and~$j$ 
and apply the antisymmetry of the Poisson bracket, to obtain
\begin{eqnarray*}
\phantom{\frac{dI_v(t)}{dt}} 
&\phantom{=}&
\sum_{i=0}^{N-1} 
\left( \sum_{j,k=0}^{N-1} H\big(v\big(t+\tfrac{j}{N}\big)\big) \, K\big(v\big(t+\tfrac{k+i}{N}\big)\big) \right)
               dK \big(v\big(t+\tfrac{i}{N}\big)\big) \, X_H\big(v\big(t+\tfrac{i}{N}\big)\big) \\
&-&\sum_{i=0}^{N-1} 
\left( \sum_{j,k=0}^{N-1} H\big(v\big(t+\tfrac{j+i}{N}\big)\big) \, K\big(v\big(t+\tfrac{k}{N}\big)\big) \right)
               dK\big(v\big(t+\tfrac{i}{N}\big)\big) \, X_H \big(v\big(t+\tfrac{i}{N}\big)\big) .
\end{eqnarray*}
This vanishes because the $1$-periodicity of~$v$ implies that for each $i$ both of the sums in the large brackets are equal to
$$
\sum_{j,k=0}^{N-1} H\big(v\big(t+\tfrac{j}{N}\big)\big) \, K\big(v\big(t+\tfrac{k}{N}\big)\big).
\phantom{\sum \sum \sum \sum}
$$
The lemma follows.

\section{Proof of Theorem~\ref{t:floer}} 
\label{s:floer}

We focus on explaining the additional ingredients compared to Floer's proof of the un-delayed case.
Following Floer~\cite{Flo89:hol} we choose an $\omega$-compatible almost complex structure~$J$ on~$M$
(that is, $\omega (v,Jw)$ is a Riemannian metric).
On the component~$\L_{\contr}$ of contractible free loops in~$M$
consider the $L^2$-inner product 
$$
g_v (\hat v_1, \hat v_2) \,:=\, \int_0^1 \omega_v \bigl( \hat v_1(t), J\2 \hat v_2(t) \bigr)\, dt 
$$
where $\hat v_1, \hat v_2$ are smooth vector fields along the loop~$v \in \L_{\contr}$.
The computation of the differential of~$\A$ before Lemma~\ref{le:crit} applies verbatim to time-dependent Hamiltonians, and yields
$$
d \A(v) \,\hat v \,=\, \int_0^1 \om \bigl( \hat v(t), \dot v(t) - \ast_v \bigr) dt 
$$
where 
$$
\ast_v \,=\, H \bigl( v(t+\tau),t+\tau \big) \cdot X_K (v(t),t) 
                     + K \bigl( v(t-\tau) ,t-\tau \bigr) \cdot X_H (v(t),t) .
$$
By the definition of the gradient and using $J^2 = - \id$ we obtain
$$
\nabla \A (v) = - J \left(\dot v(t) - \ast_v \right) .
$$
The Floer equation for $u \colon \R \times S^1 \to M$ is therefore
\begin{eqnarray} \label{e:Floer}
0 &=& \partial_s u(s,t) + \nabla \A (u(s,t)) \notag \\
&=&
\partial_s u(s,t) - J (u(s,t)) \bigl( \partial_t u(s,t) - \ast_u \bigr)
\end{eqnarray}
where
$$
\ast_u \,=\, H \bigl( u(s,t+\tau),t+\tau \big) \cdot X_K (u(s,t),t) 
                     + K \bigl( u(s,t-\tau) ,t-\tau \bigr) \cdot X_H (u(s,t),t) .
$$
Note that \eqref{e:Floer} is a \boxed{\text{\rm \bf non-local}} perturbed Cauchy--Riemann equation
in which, crucially, the non-locality is only in the lower order term~$\ast_u$.
Further, $M$ being compact, the term $\ast_u$ is uniformly bounded.
Therefore, the usual bubbling-off analysis and our assumption that $(M,\omega)$ is symplectically aspherical
imply that the space of solutions to~\eqref{e:Floer} is compact up to $s$-shift and breaking.
Since the Hessian of~$\A$ has trivial kernel, the moduli space of unparametrized broken
trajectories can be given the structure of the zero locus of a Fredholm section~$\sigma$ 
in an $M$-polyfold bundle. 
By Theorem~5.25 of~\cite{HWZ09}, the section~$\sigma$ can be perturbed to a section~$\tilde \sigma$ 
that is transverse to the zero section of the $M$-polyfold bundle.
The usual arguments in Floer theory now show that counting $0$-dimensional components 
of~$\tilde \sigma^{-1}(0)$ defines a boundary operator on the $\Z_2$-vector space generated by
the critical points of~$\A$.
Let $\HF (H,K,\tau)$ be the resulting homology.
Applying Floer's continuation argument to a homotopy from $K(x,t-\tau)$ to the constant function~$1$
(and, if necessary, from $H$ to a Floer-regular~$\widetilde H$),
we find that $\HF (H,K,\tau)$ is isomorphic to~$\HF (H)$,
and therefore, by Floer's theorem from~\cite{Flo88:Lag, Flo89:hol}, isomorphic to $\H (M;\Z_2)$.
Theorem~\ref{t:floer} now follows from the Morse inequalities.

\begin{Rmk}
{\rm
By the same reasoning, the Arnold conjecture can be proven also for the
more general delay equations discussed in this paper, 
and for any Hamiltonian delay equation in which
the delay~$\tau$ is constant or just depends on time, 
but not on the loop~$v$.
}
\end{Rmk}


%
\end{document}